\documentclass[11pt]{article}
\usepackage{epsfig} 
\usepackage{float}
\usepackage{subfigure}
\usepackage{graphicx}
\usepackage{amsmath}
\usepackage{amssymb}
\begin{document}
\title{Diffusion on Fractal Ces\`{a}ro Curve}
\author{   Alireza K. Golmankhaneh $^{a*}$}
\date{}
\maketitle \vspace{-9mm}
\begin{center}
$^a$ Young Researchers and Elite Club, Urmia Branch, Islamic Azad University, Urmia, Iran\\
 \emph{\dag E-mail address}: alirezakhalili2002@yahoo.co.in
 \end{center}

\begin{abstract}
In this paper, we apply $F^{\alpha}$-calculus on fractal Koch and Ces\`{a}ro curves with different dimensions. Generalized  Newton's second law on the fractal Koch and Ces\`{a}ro curves is suggested. Density of moving particles which absorbed  on fractal Ces\`{a}ro are derived. More,  illustrative examples are given to present the details of $F^{\alpha}$-integrals and $F^{\alpha}$-derivatives.

\vspace{.5cm} \noindent {\it {Keywords:}} $F^{\alpha}$-calculus;  Fractal Koch curve; Staircase function; Fractal Ces\`{a}ro curve
\end{abstract}
\section{Introduction}
Fractional calculus which include derivatives and integrals  with arbitrary orders is applied in science, engineering  and etc. \cite{Pod,Uchaikin-1,Balea-1}. The fractional derivatives are used to model non-conservative systems, processes with the memory effect and anomalous diffusion   \cite{Golmankhaneh-1,Golmankhaneh-3,Herrmann-2}. The fractional derivatives are  non-local but the most measurements in physics are local \cite{Hilfer-2}. As a result, in view of the fractional local derivatives and Chapman-Kolmogorov condition new Fokker-Planck equation is given   \cite{Kolwankar-1}. The local fractional derivatives lead to new measure on fractal sets \cite{Kolwankar-758}. Fractal geometry which is the generalized Euclidean geometry  has important role in science, engineering and medical science. Fractals are the shapes which have self-similar properties and fractional dimensions. Many methods have been used to build analysis on fractal sets and processes \cite{Nonnenmacher,Mandelbrot-1,Kigami-1,Strichartz,Falconer-1,Balankin-106-k,Balankin-10ll6-k}.

Recently, $F^{\alpha}$-calculus is suggested in the seminal paper as a framework by A. Parvate and A. D. Gangal which is the generalized standard calculus. $F^{\alpha}$-calculus is  calculus on fractals with the algorithmic property \cite{Gangal-1,Gangal-2,Gangal-3}. Researchers have explored this area giving new insight into $F^{\alpha}$-calculus \cite{Golmankhaneh-1015-k,Golmankhaneh-106-k,Golmankhaneh-105589-k,Golmankhaneh-105-k}. The fractal Cantor sets are considered as grating in the diffraction phenomena \cite{Golmankhaneh-104-k}. Regarding the above mentioned research we apply the   $F^{\alpha}$-calculus on fractal curves in the case of fractal Koch curves. Differential equation corresponding for a  motion of the particle on  fractal  curves is purposed. \\
Outline of the paper is as follows:\\
In  Section \ref{int-1} we summarize $F^{\alpha}$-calculus on fractal Koch and Ces\`{a}ro curves without proofs. In Section \ref{int-8}, we give our new result in this paper which includes equation of motion of the particles. Section \ref{int-8} contains our conclusion.

\section{Preliminaries \label{int-1}}
\begin{figure}[H]
    \centering
    \subfigure[ $f$ on fractal Koch curve ]
    {
        \includegraphics[scale=0.3]{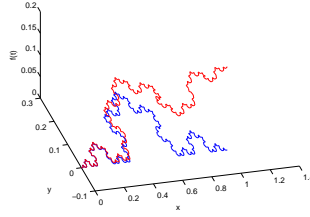}
    }
    \\
    \subfigure[ $\int f(t)d_{F}^{\alpha}t$ on fractal Koch curve ]
    {
        \includegraphics[scale=0.3]{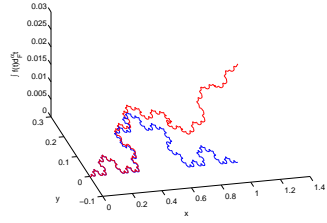}
    }
    \subfigure[ $D_{F}^{\alpha}f(t)$ on fractal Koch curve ]
    {\includegraphics[scale=0.3]{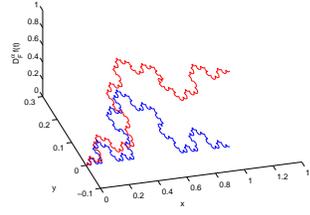}
        }
    \caption{We plot  $F^{\alpha}$-integral and $F^{\alpha}$-derivative of $f$  on fractal Koch curve with dimension 1.26 }
    \label{fig:6211}
\end{figure}

In this section, we summarize   $F^{\alpha}$-calculus on parameterize fractal curves and use in the case of fractal Koch and Ces\`{a}ro curves (see for review Refs. \cite{Gangal-3,Golmankhaneh-103-k}.\\
\textbf{Calculus on fractal Koch and Ces\`{a}ro curves:}\\ Let us consider fractal Koch and Ces\`{a}ro curves which are denoted by $F\subset \Re^3 $ and define corresponding staircase function. Fractal Koch and Ces\`{a}ro  curves are called continuously parameterizable  if there exists a function $\textbf{w}(t): [a_{0},b_{0}]\rightarrow F,~a_{0}, b_{0}\in \Re $ which is continuous one to one and onto $F$ \cite{Gangal-3,Golmankhaneh-103-k}.\\
\textbf{Definition :}~For the fractal  curves $F$ and a
subdivision $P_{[a,b]}, a<b, [a,b] \subset [a_{0},b_{0}]$ mass function is defined \cite{Gangal-3}
\begin{equation}\label{23}
  \gamma^{\alpha}(F,a,b)=\lim_{\delta\rightarrow 0} \inf_{\{P_{[a,b]}:|P|\leq
    \delta\}}\sum_{i=0}^{n-1}\frac{|\textbf{w}
    (t_{i+1})-\textbf{w}(t_{i})|^{\alpha}}{\Gamma(\alpha+1)},
\end{equation}
where $|.|$ indicates the Euclidean norm on $R^{3}$ and $P_{[a,b]}=\{a=t_{0},...,t_{n}=b\}.$ \\
\textbf{Definition}:~The staircase functions for fractal Koch and Ces\`{a}ro curves are defined
\begin{equation}\label{yyyby}
S_{F}^{\alpha}(t)=\begin{cases}
\gamma^{\alpha}(F,p_{0},t)~~~t\geq p_{0},\\
-\gamma^{\alpha}(F,t,p_{0})~~~t< p_{0},
\end{cases}
\end{equation}
where $p_{0}\in [a_{0},b_{0}]$ is arbitrary point.\\
\begin{figure}
    \centering
    \subfigure[ Fractal Koch curve and $S_{F}^{1.26}(t)$]
    {
        \includegraphics[scale=0.3]{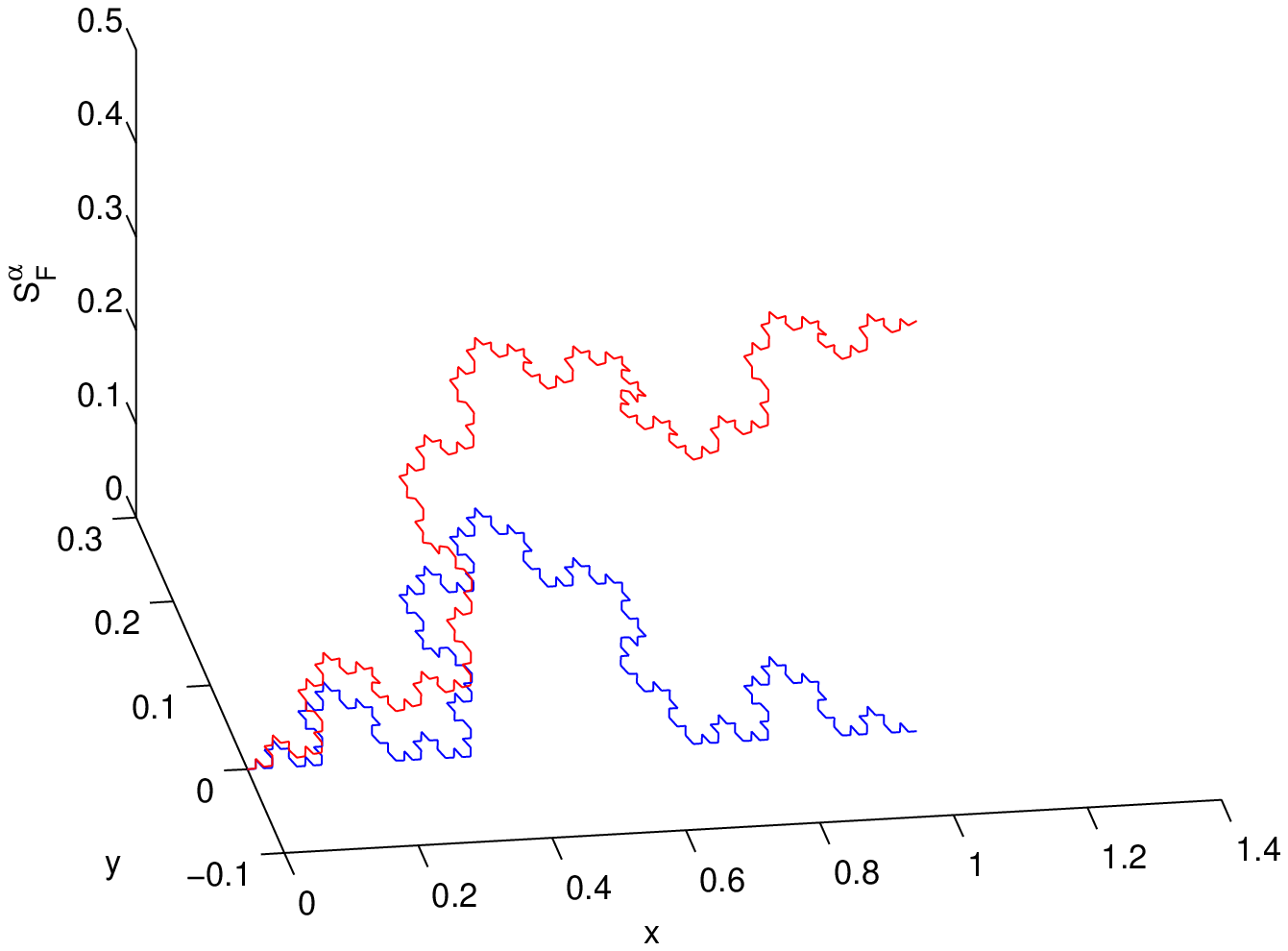} 
    }
    \subfigure[Fractal Ces\'{a}ro curve and $S_{F}^{1.7}(t)$ ]
    {
        \includegraphics[scale=0.3]{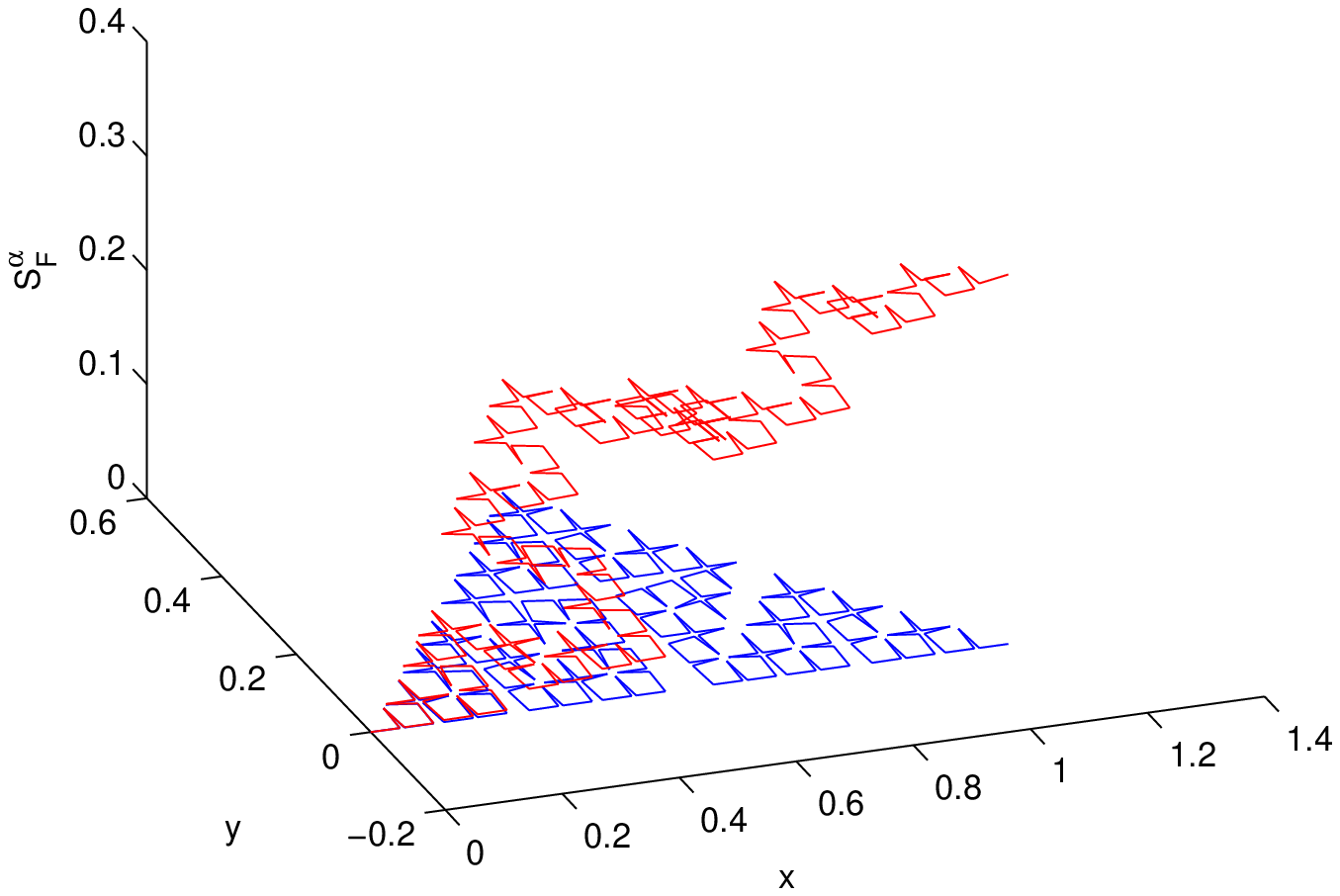} 
    }
    \caption{We have plotted fractal Koch curve and fractal Ces\'{a}ro curve in blue and corresponding staircase functions in red.}
    \label{fig:807c}
\end{figure}
In Figure [\ref{fig:807c}] we have sketched   fractal Koch and Ces\'{a}ro curves  and $S_{F}^{\alpha}(t)$ setting $\alpha=1.26$ and $\alpha=1.78$.\\
\textbf{Definition:}~The $\gamma$-dimension of fractal Koch and Ces\'{a}ro curves ($F$) are defined
\begin{align}
&\dim_{\gamma}(F)=\inf\{\alpha:\gamma^{\alpha}(F,a,b)=0\}\nonumber\\&=\sup
\{\alpha:~\gamma^{\alpha}(F,a,b)=\infty\}.
\end{align}
\textbf{Definition:~}$F^{\alpha}$-derivative of function $f$ at $\theta\in F$  is defined
\begin{equation}\label{sea}
    D_{F}^{\alpha}f(\theta)=F-\lim_{\theta'\rightarrow\theta}
    \frac{f(\theta')-f(\theta)}{J(\theta')-J(\theta)},
\end{equation}
where $J(\theta)=S_{F}^{\alpha}(w^{-1}(\theta)), \theta\in F$ and
if the limit exists \cite{Gangal-3,Golmankhaneh-103-k}.\\
\textbf{Definition:~} A number $l$ is $F$-limit of the function $f$  if we have
\begin{equation}\label{m}
    \theta'\in F ~~\textmd{and}~~~
    |\theta'-\theta|<\delta\Rightarrow|f(\theta')-l|<\epsilon.
\end{equation}
If such a number exists \cite{Gangal-3}. It is indicate by
\begin{equation}\label{z}
    l=F-\lim_{\theta'\rightarrow\theta} f(\theta').
\end{equation}
 A segment $C(t_{1},t_{2})$ of fractal Koch and  Ces\'{a}ro curve is define as
\begin{equation}\label{b}
    C(t_{1},t_{2})=\{w(t'):t'\in [t_{1},t_{2}]\},
\end{equation}
and $M,~m$ are defined as follows \cite{Gangal-3,Golmankhaneh-103-k}
\begin{equation}\label{qa}
    M[f,C(t_{1},t_{2})]=\sup_{\theta\in C(t_{1},t_{2})} f(\theta),
\end{equation}
\begin{equation}\label{x}
    m[f,C(t_{1},t_{2})]=\inf_{\theta\in C(t_{1},t_{2})}f(\theta).
\end{equation}
\textbf{Definition:~ }The upper and the lower $F^{\alpha}$-sum for the function $f$ over the subdivision $P$ are defined
\begin{equation}\label{z}
    U^{\alpha}[f,F,P]=\sum_{i=0}^{n-1}M[f,C(t_{i},t_{i+1})][S_{F}^{\alpha}(t_{i+1})
    -S_{F}^{\alpha}(t_{i})],
\end{equation}
\begin{equation}\label{awz}
    L^{\alpha}[f,F,P]=\sum_{i=0}^{n-1}m[f,C(t_{i},t_{i+1})][S_{F}^{\alpha}(t_{i+1})
    -S_{F}^{\alpha}(t_{i})].
\end{equation}
\textbf{Definition:~}$F^{\alpha}$-integral of the function $f$ is defined
\begin{align}
   &\int_{C(a,b)}f(\theta)d_{F}^{\alpha}\theta =  \underline{\int_{C(a,b)}}f(\theta)d_{F}^{\alpha}\theta=\sup_{P_{[a,b]}}L^{\alpha}[f,F,P]
  \nonumber\\& = \overline{\int_{C(a,b)}}f(\theta)d_{F}^{\alpha}\theta=\inf_{P_{[a,b]}}U^{\alpha}[f,F,P].
\end{align}
\textbf{Fundamental theorems of $F^{\alpha}$-calculus}:\\
\textbf{First Part:} If $f:F\rightarrow R$ is $F^{\alpha}$-differentiable function and
$h:F \rightarrow R$ is $F$-continuous such that
$h(\theta)=D_{F}^{\alpha}f(\theta)$, then we have
\begin{equation}\label{4po}
    \int_{C(a,b)}h(\theta)d_{F}^{\alpha}\theta=f(w(b))-f(w(a)).
\end{equation}
\textbf{Second part:} If $f$ is bounded, $F$-continuous on $C(a,b)$ and
$g: F\rightarrow R$ then we have
\begin{equation}\label{ser}
    g(w(t))=\int_{C(a,t)}f(\theta)d_{F}^{\alpha}\theta,~~~~t\in[a,b],
\end{equation}
 where we have
\begin{equation}\label{sertqw}
    D_{F}^{\alpha}g(\theta)=f(\theta).
\end{equation}
For the proofs we refer the reader to \cite{Gangal-3}.\\
\textbf{Some of the properties:}\\
1) If $f(\theta)=k\in R$ then we have $D_{F}^{\alpha}f=0$ .\\
2) If $f$ is a $F$-continuous and $D_{F}^{\alpha}f=0$ then $f=k$.\\
3) Generalized Taylor Series on fractal Koch curves is
\begin{equation}\label{mk}
    h(\theta)=\sum_{n=0}^{\infty}\frac{(J(\theta)-J(\theta'))^{n}}
    {n!}(D_{F}^{\alpha})^{n}h(\theta'),~~~~\theta\in F.
\end{equation}
\begin{align}\label{aaa}
   & \int_{C(a,b)}f(\theta)d_{F}^{\alpha}\theta=
    \int_{C(a,b)}1d_{F}^{\alpha}\theta\nonumber\\&=S_{F}^{\alpha}(b)-S_{F}^{\alpha}(a)=J((w(b))-J((w(a)),
\end{align}
where $f(\theta)=1$ is constant function \cite{Gangal-3,Golmankhaneh-103-k}.\\
\textbf{Note:}~$F^{\alpha}$-derivative and $F^{\alpha}$-integral on fractal Koch curves are  linear operators.\\
\textbf{Example 1.} Consider $f(t): F \rightarrow R$ on the fractal Koch curves as
\begin{equation}\label{b}
 f(t)=(S_{F}^{\alpha}(t))^2.
\end{equation}
The $F^{\alpha}$-derivative and the $F^{\alpha}$-integral of $f$ are
\[
\int_{C(0,t)}f(t)d_{F}^{\alpha}t=\frac{(S_{F}^{\alpha}(t))^3}{3}+k,
\]
and
\[
D_{F}^{\alpha}f(t)=2~S_{F}^{\alpha}(t),
\]
where $k$ is constant. Figure [\ref{fig:6211}] shows the graphs of  $f$, $F^{\alpha}$-integral of $f$, and $F^{\alpha}$-derivative of $f$.

\section{Equation of  motion on  fractal curves \label{int-8}}
~Generalized  Newton's second law on fractal Koch and Ces\'{a}ro curves is suggested
\begin{equation}\label{vc}
  m(D_{F}^{\alpha})^2\textbf{r}_{F}^{\alpha}(t)=\textbf{f}^{\alpha}_{F},
\end{equation}
where $\textbf{r}_{F}^{\alpha}:F\rightarrow\Re$,~ $\textbf{v}_{F}^{\alpha}(t)=D_{F}^{\alpha}~\textbf{r}_{F}^{\alpha}$ and $\textbf{a}_{F}^{\alpha}(t)=(D_{F}^{\alpha})^2~\textbf{r}_{F}^{\alpha}$ are called  generalized position,  generalized velocity and generalized acceleration  on fractal Koch and Ces\'{a}ro curves, respectively.\\
\textbf{Example 2.} Consider a force $\textbf{f}^{\alpha}_{F}=k (\hat{i}+ \hat{j})~[~M L^{\alpha}T^2~]$ such that $\textbf{f}^{\alpha}_{F}:F\rightarrow \Re$ apply on a particle with mass $m$ on fractal Koch curves.  One sees immediately that  generalized acceleration, velocity and position are
\begin{eqnarray}
\textbf{a}_{F}^{\alpha}(t)&=&\frac{k}{m},\nonumber\\~~~ \textbf{v}_{F}^{\alpha}(t)&=&\frac{k}{m}S_{F}^{\alpha}(t)+\textbf{v}_{F}^{\alpha}(0), \nonumber\\ ~~\textbf{r}_{F}^{\alpha}(t)&=&\frac{k}{2m}S_{F}^{\alpha}(t)^2+\textbf{v}_{F}^{\alpha}(0)S_{F}^{\alpha}(t)+\textbf{r}_{F}^{\alpha}(0).\label{fvr}
\end{eqnarray}
\begin{figure}[H]
  \centering
  \includegraphics[scale=0.3]{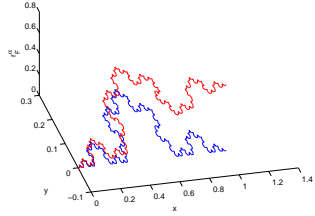}
  \caption{We have sketched $\textbf{r}_{F}^{\alpha}(t)$ for   $\textbf{v}_{F}^{\alpha}(0)=\frac{k}{2m}=1,~\textbf{r}_{F}^{\alpha}(0)=0$  }\label{thgr}
\end{figure}
We preset the graph of the Eq. (\ref{fvr}) in Figure [\ref{thgr}]

\textbf{Example 3.} Consider particles moving along  the fractal Ces\`{a}ro curve  which absorb the particles. The mathematical model for this phenomenon is given by
\begin{equation}\label{b9}
  D_{F}^{\beta}\zeta(t)=-k \zeta(t),~~~\beta=1.78.
\end{equation}
where $\zeta$ is the density of particles on fractal Ces\`{a}ro curve. Using $F^{\alpha}$-integral, it is easy to obtain the solution
\begin{equation}\label{xd}
  \zeta(t)=\zeta(0)e^{-k S_{F}^{\beta}(t)}.
\end{equation}
\begin{figure}[H]
  \centering
  \includegraphics[scale=0.3]{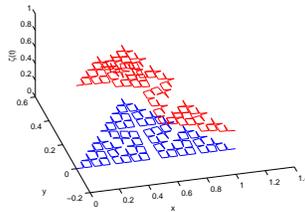}
  \caption{We plot  the density of particles for the flux of particles on fractal Ces\`{a}ro curve with absorbtion }\label{tr}
\end{figure}
Figure [\ref{tr}]  shows the graph of $\zeta(t)$ on fractal Ces\`{a}ro curve.

\section{Conclusion \label{int-7}}
In this paper,  $F^{\alpha}$-calculus  is the generalization of the standard calculus on the fractals with fractional dimension and self-similar properties. In the sense of the standard calculus the fractal Koch and Ces\`{a}ro curves are not differentiable and integrable. The $F^{\alpha}$-calculus is used to define $F^{\alpha}$-integral  and $F^{\alpha}$-derivative  on fractal Koch and Ces\`{a}ro curves. Some illustrative examples are given for presenting the details. Finally, generalized differential equation corresponding to the motions on the fractal Koch and Ces\`{a}ro curves are suggested and solved.


\begin{thebibliography}{99}
\bibitem{Pod}I. Podlubny, \textit{Fractional differential equations}, Academic Press, New York, 1999.
\bibitem{Uchaikin-1}V. V. Uchaikin,  \textit{Fractional derivatives for physicists and engineers,} Springer, Berlin, 2013.
\bibitem{Balea-1}D. Baleanu, K. Diethelm, E. Scalas, J.J. Trujillo,  \textit{Models and numerical methods}, World Scientific, Berlin, 2012.
\bibitem{Golmankhaneh-1} A.K.  Golmankhaneh, Investigations in Dynamics: With Focus on Fractional Dynamics, Lap Lambert, Academic Publishing, Germany, 2012.
    \bibitem{Golmankhaneh-3} D. Baleanu, A.K.  Golmankhaneh,  A.K.  Golmankhaneh, R.R.  Nigmatullin,  Newtonian law with memory,  \textit{Nonlinear Dyn.,}  60(1-2), (2010), 81-86.
\bibitem{Herrmann-2} R. Herrmann, \textit{Fractional calculus: an introduction for physicists,} World Scientific, 2014.
\bibitem{Hilfer-2}  R. Hilfer,  ed., \textit{Applications of fractional calculus in physics}, World Scientific, 2000.
\bibitem{Kolwankar-1} K. M. Kolwankar, A.D. Gangal, Local fractional Fokker-Planck equation,\textit{ Phys. Rev. Lett.} 80(2) (1998), 214-217.
\bibitem{Kolwankar-758} K. Kolwankar, J. L. Vehel, \textit{Measuring functions smoothness with local fractional derivatives}, Fractional Calculus and Applied Analysis, 4(3) (2001), 285-301.
\bibitem{Nonnenmacher} F. T. Nonnenmacher, A. L.  Gabriele, R. W.  Ewald, eds. \textit{Fractals in biology and medicine}. Birkhäuser, 2013.
\bibitem{Mandelbrot-1} B.B.  Mandelbrot, \textit{The Fractal Geometry of Nature}, Freeman and Company, 1977.
\bibitem{Kigami-1} J.  Kigami, \textit{Analysis on fractals} , Volume 143 of Cambridge Tracts in Mathematics, Cambridge University Press, Cambridge, 2001.
\bibitem{Strichartz} R. S.  Strichartz, \textit{Differential equations on fractals: a tutorial}, Princeton University Press, 2006.
\bibitem{Falconer-1} K.  Falconer,   \textit{Techniques in Fractal Geometry}, John Wiley and Sons, 1997.
\bibitem{Balankin-106-k} A. S. Balankin, Effective degrees of freedom of a random walk on a fractal, \textit{Phys. Rev. E,} 92(6) (2015), 062146.
\bibitem{Balankin-10ll6-k} F. Alonso-Marroquin, P. Huang, D. A. Hanaor, E. A. Flores-Johnson, G. Proust, Y. Gan, and L. Shen.,Static friction between rigid fractal surfaces, \textit{Phys. Rev. E} 92(3) (2015), 032405.
\bibitem{Gangal-1} A. Parvate, A.D. Gangal,  Calculus on fractal subsets of real-line I: Formulation,
     \textit{Fractals,} 17(1), (2009), 53-148.
\bibitem{Gangal-2}A. Parvate, A.D. Gangal,  Calculus on fractal subsets of real-line II: Conjugacy with ordinary calculus, \textit{Fractals,} 19(3) (2001), 271-290.
\bibitem{Gangal-3} A. Parvate, S. Satin,  A.D. Gangal, Calculus on fractal curves in $R^n$,  \textit{Fractals} 19(1) (2011), 15-27.
\bibitem{Golmankhaneh-1015-k} A.K. Golmankhaneh, V. Fazlollahi, D. Baleanu,  Newtonian mechanics on fractals subset of real-line, \textit{Rom. Rep. Phys,} 65, (2013) 84-93.
\bibitem{Golmankhaneh-106-k} A.K. Golmankhaneh, D. Baleanu, Non-local Integrals and Derivatives on Fractal Sets with Applications, \textit{Open Phys.}, 14(1) (2016), 542-548.
    \bibitem{Golmankhaneh-105589-k} A. K. Golmankhaneh, C. Tunc, On the Lipschitz condition in the fractal calculus, \textit{Chaos, Soliton  Fract.}, 95 (2017), 140-147.
 \bibitem{Golmankhaneh-105-k} A.K. Golmankhaneh,  D. Baleanu, Fractal calculus involving Gauge function, \textit{Commun. Nonlinear Sci.}, 37 (2016), 125-130.
 \bibitem{Golmankhaneh-104-k} A. K.  Golmankhaneh, D. Baleanu, Diffraction from fractal grating Cantor sets, \textit{J. Mod. Optic,} 63(14) (2016), 1364-1369.
\bibitem{Golmankhaneh-103-k} A.K. Golmankhaneh,  A.K. Golmankhaneh, D. Baleanu, About Schr\"{o}odinger equation on fractals curves imbedding in $R^3$ , \textit{Int. J.  Theor.  Phys.,} 54(4) (2015), 1275-1282.
\end{thebibliography}

\end{document}